\documentclass{amsart}
\usepackage{amssymb, graphicx}
\newcommand{\ncmnd}{\newcommand}
\ncmnd{\nthm}{\newtheorem}
\ncmnd{\R}{\mathbf{R}}

\DeclareMathOperator{\inj}{inj}
\DeclareMathOperator{\Conv}{Conv}

\theoremstyle{plain}

\nthm{thm}{Theorem}
\nthm{cor}[thm]{Corollary}
\nthm{prop}[thm]{Proposition}
\nthm{lem}[thm]{Lemma}
\nthm{examp}[thm]{Example}
\theoremstyle{definition}
\nthm{note}[thm]{Note}
\nthm{rem}[thm]{Remark}
\nthm{defn}[thm]{Definition}

\newcounter{step}

\begin{document}

\title{Gauss Equation and Injectivity  Radii for Subspaces\\in Spaces of
Curvature Bounded Above}
\author{Stephanie B. Alexander}
\address{Department of Mathematics, University of Illinois\\1409 W. Green St., Urbana, Illinois 61801}
\email{sba@math.uiuc.edu}
\author{Richard L. Bishop}
\address{Department of Mathematics, University of Illinois\\1409 W. Green St., Urbana, Illinois 61801}
\email{bishop@math.uiuc.edu}

\begin{abstract}

A Gauss Equation is proved for subspaces of Alexandrov spaces of curvature bounded above by $K$.  That is, a subspace of extrinsic curvature $\le A$, defined by a cubic inequality on the difference of arc and chord, has intrinsic curvature $\le K+A^2$.  Sharp bounds on injectivity radii of subspaces, new even in the Riemannian case, are derived.

\end{abstract}

\maketitle

\section{Introduction}
\label{sec:intro}

Alexandrov spaces are metric spaces with curvature bounds in
the sense of local triangle comparisons with constant curvature spaces.
In this paper, we consider  spaces of
curvature bounded above (CBA), and their global counterparts,
CAT($K$) spaces.  Examples include Riemannian manifolds with upper
sectional curvature bounds,  possibly with boundary (\cite{ABB1}), polyhedra with link conditions, and Tits boundaries and 
asymptotic cones (see \cite{BH}).  A
key property of CAT($K$)
spaces is their preservation under Gromov-Hausdorff convergence. 
CAT($K$)  spaces
are appropriate target spaces in harmonic
map theory (see, for example, \cite{EF,GS,J}), and play an
important role in geometric group theory (see \cite{BH}).

Analogues of the Gauss Equation, governing the passage of curvature
bounds to subspaces from ambient spaces, tend to be challenging in
Alexandrov spaces.  For instance,  a major unsolved problem in the
theory of spaces of curvature bounded below  is whether the boundary of a
convex set inherits the curvature bound.  Somewhat more is known  for
CBA.  A classical theorem of Alexandrov states that a curvature bound above
is inherited by ruled surfaces \cite{A}.   It is an open problem whether
saddle surfaces inherit such a bound, but Mese showed that minimal surface
immersions do so \cite{M},  and Petrunin \cite{P}, that
``metric minimizing'' surfaces do so.

Recently Lytchak proved that if $M$ has curvature bounded above by $K$, and
$N$ is a complete subset for which there exists $\rho>0$ such that
intrinsic distances $d_N = s$
and extrinsic distances $d_M = r$ satisfy  $s-r\le Cr^3$ for $r<\rho$ (i.e., 
$N$ is  \emph{(C,2,$\rho$)-convex}), 
then $N$ has \emph{some} intrinsic curvature
bound above \cite{L1}. For a subset $N$ of a Riemannian manifold,
an equivalent condition on $N$ is \emph{positive 
reach}, namely, uniform neighborhoods in which $N$ has the unique footpoint 
property \cite{L2}.  In general, let us say \emph{$N$ is a subspace of 
extrinsic curvature $\le A$} in $M$ if there is a length-preserving map 
$F:N\to M$ between intrinsic metric spaces, where  $N$ is complete and
\begin{equation}\label{eq:ext}
s-r \le\frac{A^2}{24} r^3+ o(r^3) 
\end{equation}
for all pairs of points having $s$ sufficiently small. (For Riemannian submanifolds,  this  is equivalent to a bound, $|II| \le A$, on the second fundamental form.) Then points of $N$ have $(C,2,\rho)$-convex neighborhoods.  It follows from \cite{AB3} that points of $N$ have neighborhoods in which $r$ is at least the chordlength of an arc of constant curvature $A$ and length $s$ in the model plane $S_K$ (see Theorem \ref{thm:CBC1}.1 below).

In this paper, we extend the Gauss Equation to Alexandrov spaces of
curvature bounded above, by proving the following sharp bound for 
subspaces of extrinsic curvature $\le A$.

\begin{thm}[Gauss Equation]
          \label{thm:Gausseq} 
Suppose $N$ is a
subspace of extrinsic
curvature $\le A$ in an Alexandrov space of curvature bounded above
by $K$. Then $N$ is
an Alexandrov space of curvature bounded above by $K+A^2.$
\end{thm}

\begin{rem} This bound is realized by  constantly curved hypersurfaces of Euclidean,
spherical and hyperbolic spaces.  At first one might think that Riemannian submanifolds of higher codimension offer a counterexample to this theorem, and that the correct bound should be $K+2A^2$.  On closer inspection, however, one sees that for any plane section, normals to the submanifold may be chosen so that at most two of the corresponding subdeterminants of $II$ are nonzero and one of them is nonpositive.  Therefore for Riemannian submanifolds, while the sharp \emph{lower} bound is $K-2A^2$ when ambient curvature is $\ge K$, the sharp \emph{upper} bound is $K+A^2$ when ambient curvature is $\le K$.
\end{rem}

There are important classes of subspaces for which we can compute
sharp extrinsic curvature bounds, and hence sharp intrinsic curvature
bounds by Theorem \ref{thm:Gausseq}. Fibers of warped products are
such a class.  Warped products of Alexandrov spaces extend standard
cone and suspension constructions from one-dimensional to arbitrary
base, and gluing constructions from $0$-dimensional to arbitrary fiber
\cite {AB1}; we expect them to be a major source of constructions and
counter-examples in the Alexandrov setting.  Theorem \ref{thm:Gausseq}
allows us to calculate the intrinsic curvature bound of the fiber of a
CAT($K$) warped product, as we shall discuss in detail elsewhere.

Another significant application of Theorem \ref{thm:Gausseq} is to 
injectivity radii. Theorem \ref{thm:injrad} gives  a sharp estimate on
the injectivity radius of a
subspace of bounded extrinsic curvature, in  terms of the circumference
$c(A,K)$ of a circle of curvature $A$ in the simply connected,
$2$-dimensional space form  $S_K$ of curvature $K$.
As always, we set $\pi/k =\infty$ if $k\le 0$.

\begin{thm}\label{thm:injrad} 
 Suppose $N$ is a subspace of extrinsic curvature $\le A$ in a CAT($K$) space.  Then
\begin{equation}\label{eq:inj}
\inj _N \ge \min\{
\frac{\pi}{\sqrt{K + A^2}}, \frac{1}{2} c(A,K)\}.
\end{equation}
\end{thm}

The definition of the injectivity radius $\inj _N$ agrees with the usual one when $N$ is Riemannian or locally compact.  In the non-locally compact case, a slightly stronger definition is more appropriate (see \S \ref{sec:injbd}). Even in the case of Riemannian manifolds, our estimates on the
injectivity radius of a submanifold are new as far as we know.  Much
weaker dimension-dependent estimates have been used in \cite{Co} and
\cite{S}. The existence of some dimension-independent bound in the
general case is proved in \cite{L2}.

The following corollary holds, in particular, for Riemannian
submanifolds with $|II|\le A$ in a Hadamard manifold. Part 2 is an
``immersion implies embedding'' theorem, which in the case of
Riemannian hypersurfaces appeared in \cite{Ar}.

\begin{cor} \label{cor:CAT0}
Let $N$ be a subspace of extrinsic curvature $\le A$
in a CAT(0) space $M$.
\begin{list}
{\emph{\arabic{step}.}}
{\usecounter{step}
\setlength{\rightmargin}{2\leftmargin}}
\item
$N$ has injectivity radius at least
$\pi/A$, and any closed ball of radius $\pi/2A$ in $N$ is CAT($A^2$).

\item
If $M$ is CAT($-A^2$), then $N$  is CAT(0) and embedded.  
\end{list}
\end{cor}

An interesting example to which Theorems \ref{thm:Gausseq} and
\ref{thm:injrad} apply is that of tubular neighborhoods of convex
sets. Their extrinsic curvature is analyzed in \S \ref{sec:sublevels} below.  A subset is
\emph{$\pi$-totally-convex} if it contains every geodesic of length
$<\pi$ joining pairs of its points.

\begin{examp} \label{examp:tub} 
Let $T$ be a $\pi$-totally-convex set in a CAT(1) space $M$,
and $N$ be the set of points at distance $\le \rho$ from $T$,
where $\rho < \pi/2$. Then $N$ has CBA by $ \sec\rho$,
 and injectivity radius $\ge \pi\cos\rho$.
\end{examp}


\section{Outline of paper}
\label{sec:outline}

A Riemannian manifold has curvature bounded above by $K$ if and only
if the lengths of normal Jacobi fields satisfy $f'' \ge -Kf$ in the
barrier sense. (In this paper, a continuous function with this
property will be called, briefly, \emph{$K$-convex}.)  The same
formulation was extended to Riemannian manifolds with boundary in
\cite{ABB1}, and guides us now in the proof of Theorem
\ref{thm:Gausseq}.  However, there we made strong use of the smooth
Gauss Equation, which is one of the cornerstones of Riemannian
geometry.  In contrast, our task here is to find a new approach to
that theorem which holds in a far more general setting.

To estimate the curvature of the subset $N$, 
we may start with the knowledge that $N$ has some upper curvature bound, and 
with Lytchak's idea of
bounding the extrinsic curvature of
two nearby $N$-geodesics in their ambient ruled surface by projecting
to $N$ \cite{L1}. Obtaining a sharp bound for this
curvature (Step 2 below)
is the difficult core  of our paper.

A \emph{geodesic} is an isometric embedding of an interval unless otherwise described. Here are the first two steps of the proof of
Theorem \ref{thm:Gausseq}:

\begin{list}
{\textbf{Step \arabic{step}.}}
{\usecounter{step}
\setlength{\rightmargin}{2\leftmargin}}
\item For a fan of geodesics $\gamma$ in $N$ converging to a
base geodesic $\sigma$, the geodesic chords in $M$ connecting
$\gamma(t)$ to
$\sigma(t)$ become arbitrarily close to normal to $\sigma$ and $\gamma$,
forming a ruled surface $R_\gamma$ with curvature bounded above by $K$.

\item The extrinsic curvature of $\sigma$ and $\gamma$ in
$R_\gamma$ is bounded above by $A^2(1 +
\epsilon)d(\sigma(t),\gamma(t))/2$, where
$\epsilon \to 0$ as $\gamma \to \sigma$.
\end{list}

\begin{figure}[h]
\begin{center}
\includegraphics[width= 3.5in]{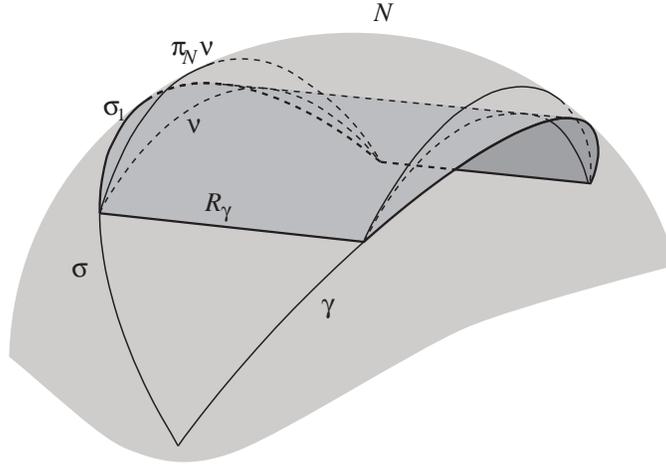}
\end{center}
\caption{Geometric layout}
\label{fig:lunes}
\end{figure}

Let us describe the geometry behind Step 2.    By definition, we must
estimate the difference in lengths between a subarc $\sigma_1$ of the
$N$-geodesic $\sigma$ and its chord $\nu$ in $R_\gamma$, where the
\emph{chord} of an arc is  the geodesic between its
endpoints. Lytchak's bound is based on
the estimate
\begin{equation} \label{eq:discrep}
\ell(\sigma_1)-\ell(\nu)\le \ell(\pi_N\nu)- \ell(\nu),
\end{equation}
where $\pi_N$ is projection to $N$. See Figure 1.  Equation
(\ref{eq:discrep}) holds
because the $N$-geodesic $\sigma_1$ is no longer than the curve
$\pi_N\nu$ with the same endpoints. 

In seeking a sharp bound,
we exploit estimates from a globalization theorem for 
arc/chord curvature in \cite{AB3} (Theorem \ref{thm:CBC1} below).  
In consequence (Lemma \ref{lem:ubproj} below), we obtain a bound for
$\ell(\pi_N\nu)-
\ell(\nu)$ which in combination with Equation (\ref{eq:discrep})
yields
an upper curvature bound for $N$ of $K+2A^2$.  In other words, it
turns out that Equation (\ref{eq:discrep}) is too weak to give the Gauss
Equation we seek; a factor of $2$ is needed on the lefthand side.
Figure \ref{fig:lunes} suggests where to  look for this factor
of $2$.   Namely, as $\ell(\sigma_1)\to 0$, it seems that
$\pi_N\nu$ may be as much longer than $\sigma_1$, as $\sigma_1$ is
longer than
$\nu$!  Indeed, this is what we prove.

We use a version for CAT($K$) spaces (Lemma \ref{lem:shrink} below)
of a standard
lemma in smooth comparison geometry.  The latter is a length
comparison between two curves with the same expression in Fermi
coordinates, in ruled surfaces generated by vector fields that are
parallel along and orthogonal to two base geodesics (see  \cite[``Corollary
of Rauch II'', p. 31-32]{CE}).  For the longer of the two
curves, the parallel condition can be dropped, since allowing the
ruled surface generators to ``twist'' as they move along the base
geodesic only lengthens a curve with given Fermi coordinates.

The remaining steps in the proof of Theorem \ref{thm:Gausseq} are:

\begin{list}
{\textbf{Step \arabic{step}.}}
{\usecounter{step}
\setlength{\rightmargin}{2\leftmargin}}
\setcounter{step}{2}
\item The distance between the respective chords in $R_\gamma$
of two subarcs of $\sigma$ and $\gamma$ is
$(K+\epsilon)$-convex.  Adding to this the width
estimates that follow from Step 2, we find that
$d(\sigma(t),\gamma(t))$ satisfies the
$(K+\epsilon+A^2(1+\epsilon))$-convexity  midpoint inequality.

\item  The  usual behavior of  normal and tangential Jacobi
field lengths extends to CBA spaces. Normal Jacobi field lengths
along $\sigma$ are $(K+A^2)$-convex
       and so, by a development argument, $N$ has CBA
by $K+A^2$.

\end{list}

In \cite{ABB1, ABB2}, Jacobi fields were studied in Riemannian
manifolds with boundary, which form a class of CBA spaces by Alexander et al.
\cite{ABB1}.  In Step 4, we must capture the notion of normal and
tangential Jacobi fields lengths in the general CBA setting.  This is
because tangential ones, which are linear, must be discarded in order
to identify negative curvature bounds.

Throughout, we may apply the following globalization theorem to $N$-geodesics, which have extrinsic curvature $\le A$ in $M$ because $N$ does.  The \emph{base
angles} of an arc $\gamma$ are the angles it makes with its chord at
their common endpoints, and the \emph{width} of $\gamma$ is the
smallest tubular radius containing $\gamma$ about the chord.  A
\emph{$k$-curve} in $S_K$ is a curve of constant extrinsic curvature
$k$; thus a complete $k$-curve is a circle, geodesic line, horocycle
or equidistant curve.

\begin{thm}[AB3]\label{thm:CBC1}
Let $\gamma$ be a curve of extrinsic curvature $\le k$ in a CAT($K$)
space.

\begin{list}
{\emph{\arabic{step}.}}
{\usecounter{step}
\setlength{\rightmargin}{\leftmargin}}

\item
If the sum of the arclength and chordlength of $\gamma$ is $<
2\pi/\sqrt{K}$, then $\gamma$ has the same arclength and chordlength
as a $k'$-curve in $S_K$, for some $k'\le k$.
\item
If $\gamma$ is closed (not necessarily closing
smoothly) and nonconstant, then $\gamma$ is no shorter than a complete
$k$-curve (necessarily a circle) in $S_K$.
\item
If $\gamma$ has length $\le$  half a  complete $k$-curve in $S_K$,
then the base angles and width of $\gamma$ are no more than they are
for a $k$-curve in $S_K$ of the same length.
\end{list}
\end{thm}

\begin{rem} The definition of \emph{extrinsic curvature $\le K$} used to 
prove the theorem above was slightly weaker than that used in this paper.
Here, inequality (\ref{eq:ext}) is imposed uniformly at all points of $N$, whereas above 
it was imposed on the distances from each point separately.  We may take the
weaker form in this paper if $N$ is a geodesic space (e.g., if $N$ is  
locally compact), since then uniformity follows from Part 1 above.
\end{rem}


\section{Majorizing Lemmas}
\label{sec:majlem}

By a \emph{CAT(K)} space, we mean here a complete metric space in
which any two points at distance $<\pi/\sqrt{K}$ are joined by a geodesic, and for any geodesic triangle
$\triangle$ of perimeter $<2\pi/\sqrt{K}$, the distances between
points of $\triangle$ are $\le$ the distances between corresponding
points on the triangle with the same sidelengths in $S_K$. A space has
\emph{curvature bounded above (CBA) by $K$} if every point has a
CAT($K$) convex neighborhood. As references see \cite{BH, BBI}. 

An important tool for studying curves in  CBA
spaces is Reshetnyak's Majorization Theorem \cite{R}:

\noindent {\bf RMT.} \emph{Let $\gamma$ be a closed curve of length
$<2\pi/\sqrt{K}$ in a CAT($K$) space $M$. Then there is a closed curve
$\widetilde\gamma$ which is the boundary of a convex region $D$ in
$S_K$ and a distance-nonincreasing map $\varphi: D \to M$ such that
the restriction of $\varphi$ to $\widetilde\gamma$ is an
arclength-preserving map onto $\gamma$.}

Such a map $\varphi$ is called a majorizing map for $\gamma$. Note
that it is an immediate consequence of the minimizing property of
geodesics that for a geodesic subarc of $\gamma$, the corresponding
subarc of $\widetilde\gamma$ is also a geodesic segment; hence, RMT is
a broad generalization of the defining property of a CAT($K$) space,
as is seen by taking $\gamma$ to be a triangle. Moreover, an
extrinsinc curvature bound at a point of $\gamma$ is inherited by the
corresponding point of $\widetilde\gamma$.

Our first application of the RMT is to prove the Fermi lemma on
curvelengths mentioned in the preceding section. The setting is a
CAT$(K)$ space $M$ where the nearest-point projection $\pi_\nu$ to a
given nontrivial geodesic segment $\nu$ is well-defined and
continuous.  (This will always hold when $K\le 0$, and holds when
$K>0$ under standard size restrictions; see \cite[p. 176-178]{BH}.)
Choose a geodesic $\widetilde\nu$ in $S_K $, where $\nu$ and
$\widetilde\nu$ are parametrized by arc length $u$ on the same
interval $[0,r],\,\, r>0$. Define $\psi: M \to S_K$ as follows:
$\psi(M)$ lies on one side of $\widetilde\nu$; if $\pi_\nu(p) =
\nu(u)$, then $d(p, \nu(u)) = d(\psi(p), \widetilde\nu(u))$, and the
geodesic segment from $\psi(p)$ to $\widetilde\nu(u)$ is normal to
$\widetilde\nu$.

\begin{lem} [Fermi Lemma]\label{lem:shrink}
The map $\psi: M \to S_K$ is distance-nonincreasing.  If
$d(\psi(p),\psi(q)) = d(p,q)$, then $p, q, \pi_\nu(q),\pi_\nu(p)$ are
the vertices of a quadrilateral region isometric to a region in $S_K$.
\end{lem}

\begin{proof}
For $p, q \in M$, let $p', q'$ be their projections to $\nu$.  If
$p'=q'$, then $\psi(p)=\psi(q)$.  If $p'\ne q'$, consider the
quadrilateral $pqq'p'$.  The angles at $p', q'$ are nonacute, so that
if we take a majorizing quadrilateral $\tilde p \tilde q \tilde q'
\tilde p'$ with base $\tilde p' \tilde q'$ at the corresponding points
of $\widetilde\nu$, in accordance with RMT, then the base angles
remain nonacute. Now we deform $\tilde p \tilde q \tilde q' \tilde p'$
so that it has the same base segment and its sides are normal to
$\widetilde\nu$. The new top vertices are $\psi(p), \psi(q)$.  The top
side is seen to be shorter than $d(p,q)$ because the deformations of
$\tilde p$ and $\tilde q$, along the circles centered at $\tilde p'$
and $ \tilde q'$ that determine the new quadrilateral, steadily
decrease the length of the top side.
\end{proof}

We also require the following local estimates on a
subset $N$ of bounded extrinsic curvature.  By               
\cite{L2},     in a sufficiently
small ambient neighborhood of a point of $N$, the nearest-point
projection $\pi_N$ is uniquely defined and continuous. 
Part 2 below sharpens Lytchak's estimate of 
$1 + Cd(x,N)$ for a Lipschitz constant for $\pi_N$ on a neighborhood
of a point $x$ \cite{L1}.

\begin{lem} \label{lem:width} 
Let $N$ be a subset of extrinsic curvature $\le k$ in a CAT$(K)$ space
$M$, and $\nu$ be a curve in $M$ parametrized by arclength $u \in
[0,r]$.

\begin{list}
{\emph{\arabic{step}.}}
{\usecounter{step}
\setlength{\rightmargin}{2\leftmargin}}

\item
If $\nu$ is a geodesic in $M$ with ends on $N$, then
$$d(\nu(u), N) \le \frac{k}{2}u(r-u) + O(r^4).$$

\item
If  $\nu$ is any curve within
the injectivity distance to $N$, then
\begin{equation}\label{eq:width1}\ell(\pi_N\nu) \le r + \int_0^r
[kd(\nu(u), N) +O(d(\nu(u), N)^2)]du.\end{equation}

\end{list}

\end{lem}

\begin{proof} \emph{1.}  Recall that an $N$-geodesic has extrinsic
curvature $\le k$ in $M$. Let $\sigma$ be an
$N$-geodesic connecting the ends of $\nu$. Consider the closed curve
formed by $\sigma$ and its chord $\nu$.  By RMT, this curve is
majorized by a map from a region in $S_K$ bounded by a convex arc of
curvature $\le k$ and its chord $\tau$ of length $r$.  Parametrizing
$\tau$ by arclength $u \in [0,r]$, we know that the distance from
$\tau(u)$ to the convex arc is at most equal to the distance from
$\tau(u)$ to the $k$-curve with the same chord $\tau$ and lying on the
same side of $\tau$ (for details, see \cite{AB3}).  Setting $v=r-u$,
the latter distance may be calculated as
$$w(u) = \frac{k}{2}uv\left(1 + \frac{k^2}{64} uv +
\frac{K}{12}(u^2 + 3uv + v^2) + O(r^4)\right).
$$ (When $u = \frac{r}{2}$ this reduces to the expression for the
width of the closed curve in terms of $r$; the corresponding formula
for width in terms of the length of the $k$-curve was given in
\cite{AB3}.)  Thus
$$
d(\nu(u), N) \le d(\nu(u), \sigma) \le w(u) = \frac{k}{2}uv + O(r^4).
$$
\end{proof}

\begin{proof} \emph{2.}  Set $d(u)= d(\nu(u), N)$.  We are going to construct
a curve in $S_K$ of length $r=\ell(\nu)$, parametrized by arclength
$u\in[0,r]$, which projects with the \emph{same} distance function
$d(u)$ into a complete $k$-curve.  Moreover, if $\widetilde\sigma$ is
the image arc of this projection, then
\begin{equation}\label{eq:proj}
\ell(\pi_N\nu)\le \ell(\widetilde\sigma).
\end{equation}

In consequence, we need only verify Equation (\ref{eq:width1}) when
$M=S_K$, $N=\widetilde\sigma$, and $\nu$ is any curve of length $r$.
In this case, the Lipschitz constants of the tangent maps of the
projection to $\widetilde\sigma$ may be calculated from the lengths of
normal Jacobi fields in $S_K$ along geodesics radiating orthogonally
to $\tilde\sigma$.  For example, if $K > 0$, we can take the length of
the Jacobi field to be $\sin(\sqrt{K}t)$; the length at points of
$\tilde\sigma$ is given by taking $t=t_1<\pi/2$, where
$\cot(\sqrt{K}t_1) = k$, while the length at distance $d$ from
$\tilde\sigma$ is given by taking $t=t_1-d$.  This gives a Lipschitz
constant of
$$\frac{\sin(\sqrt{K}t_1)}{\sin(\sqrt{K}[t_1 - d])} =
1 + kd +(k^2 +\frac{K}{2})d^2 + O(d^3).$$
The series expression for the ratio of Jacobi field lengths in terms
of $d$ and
$K$ persists for all values of $K$. The bound (\ref{eq:width1}) on 
$\ell(\pi_N\nu)$ now
follows by integrating this ratio.

To carry out the  construction, partition $\nu$ into subarcs, and
apply RMT to each closed figure whose base is a chord of that subarc, whose
sides are the minimizers from the endpoints of the base to $N$, and
whose top is the $N$-geodesic joining the footpoints of these
minimizers.  The corresponding convex curve in $S_K$ has a geodesic
base, two geodesic sides, and a top curve of curvature $\le k$ making
nonacute angles with both sides; all of these arcs have the same
length as their corresponding arcs in $M$.  Now replace the top curve
in $S_K$ by the $k$-curve with the same endpoints;  this move does
not decrease the length of the top or the angles between the top and
sides.  The latter angles are nonacute and, if not right angles, may
now be made right by hinging the sides inward, thereby reducing the
length of the base.  The final move on each figure in $S_K$ is to
lengthen the top $k$-curve, preserving its right angles with the
sides and the lengths of the sides, until the baselength is restored
to that of the corresponding segment of $\nu$.  For any partition of
$\nu$, we now glue the resulting figures in $S_K$ along corresponding
sides.  By construction, the top curves form a $k$-curve whose
length is at least that of a broken $N$-geodesic approximation of
$\pi_N(\nu)$. Thus we have constructed
a polygonal curve in $S_K$ which converges, as the original partition
of $\nu$ is refined, to a curve of the same length as
$\nu$.  By construction, this curve projects, with the same distance
function as $\pi_N|\nu$, to a $k$-curve $\widetilde\sigma$
satisfying (\ref{eq:proj}).
\end{proof}

The following lemma gives width and base angle estimates that are a
direct consequence of Theorem \ref{thm:CBC1}.3 and the power series
expansions in \cite[Remark 6.2]{AB3}:

\begin{lem}\label{lem:CBC2}
For a curve of extrinsic curvature $\le k$ in a CAT($K$) space, the
base angles $\varphi$ and width $W$ of an arc with chordlength $r$
satisfy
\begin{list}
{\emph{\arabic{step}.}}
{\usecounter{step}
\setlength{\rightmargin}{2\leftmargin}}
\item
$\varphi \le kr/2 + O(r^3),$
\item
$W \le kr^2/8 +O(r^3).$
\end{list}
\end{lem}


\section{Proof of Theorem \ref{thm:Gausseq} (Gauss Equation)}
\label{sec:Gausspf}

Now we are ready to prove our Gauss Equation, namely, a sharp
curvature bound for a subset $N$ of extrinsic curvature $\le A$ in an
Alexandrov space $M$ of CBA by $K$.  The proof breaks into Steps 1 to
4, as outlined in \S \ref{sec:outline}.


\subsection{Step 1:  The ruled surfaces $R_\gamma$}
\label{ss: pf1}

By a \emph{fan} in $N$, we mean a one-parameter family of
$N$-geodesics $\gamma=\gamma_u$ in $N$, originating at a point $p$ and
with righthand endpoints moving along a geodesic with parameter
$u$. (Recall that all geodesics are parametrized by arclength.)  Set
$\gamma_0=\sigma$.  Since $N$ has some upper curvature bound $L > 0$,
then the balls centered at $p$ are convex as long as they are
contained in a CAT($L$) neighorhood of $p$ and the distance from $p$
is $< \pi/2\sqrt{L}$; we assume that the generators $\gamma_u$ of the
fan are shorter than that radius.

The function $f(t) = d(\sigma(t), \gamma(t))$ is $L$-convex and
increasing.  Here, $L$-convexity follows from RMT and the fact that
the distance between geodesics in the model space $S_L$ is $L$-convex
for $L\ge 0$.  By extracting a subsequence $\gamma_i$ of the fan
generators, and setting $f_i=d(\sigma, \gamma_i)$, we may assume the
existence of
\begin{equation} \label{eq:jaclength}
F(t) = \lim_{i\to\infty} u_i^{-1}f_i(t).
\end{equation}
Moreover, $F$ is continuous and $L$-convex  on $[0,\ell_0)$, where $
\ell_0=|\sigma|$. (See the argument in \cite[p. 178]{ABB2}.) The
function $F$ plays the role of a normal Jacobi field \emph{length};
we do not need the notion of  Jacobi field \emph{direction}.

Recall that $L$-convex functions have one-sided derivatives, and
derivatives almost everywhere.  Since, moreover, $L$-convex functions
converge with their one-sided derivatives, it follows from Equation
(\ref{eq:jaclength}) that the one-sided derivatives $f_i'$ are $\le
Cu_i$, hence converge unformly to $0$.  By first variation, $f_i'(t_+)
= -\cos \alpha_i (t)- \cos \beta_i(t)$, where $\alpha_i (t)=
\angle(\sigma(t)\gamma_i(t)\gamma_i(\ell_0))$ and $\beta_i (t)=
\angle(\gamma_i(t)\sigma(t)\sigma(\ell_0))$. By convexity of balls,
$\alpha \ge \pi/2, \beta\ge \pi/2$. Hence $\alpha_i, \beta_i$ are
uniformly close to $\pi/2$ as $i\to\infty$. (If $f_i$ vanishes on an
initial interval, we set $\alpha_i=\beta_i=\pi/2$ there.) Similarly,
the ``supplementary'' angles $\angle(\sigma(t)\gamma_i(t)p)$ and $
\angle(\gamma_i(t)\sigma(t)p)$ are also uniformly close to
$\pi/2$. This shows that the ruled surface in $N$, formed by
connecting the pairs of points $\gamma_i(t)$ and $\sigma(t)$ by
$N$-geodesics, has its rulings nearly normal to $\gamma_i$ and
$\sigma$.

When we form the corresponding ruled surface $R_{\gamma_i}$ in $M$,
the rulings are chords of the $N$-geodesic rulings, where the latter
have extrinsic curvature $\le A$ by hypothesis. Theorem
\ref{thm:CBC1}.3 gives an upper bound for the angles between the two
kinds of rulings, for $i$ sufficiently large (depending only on $A,
K$). Therefore the rulings of $R_{\gamma_i}$ may be assumed
arbitrarily close to normal to $\gamma_i$ and $\sigma$.

By Alexandrov's theorem on ruled surfaces, $R_{\gamma_i}$ has
curvature bounded above by $K$.


\subsection{ Step 2: Extrinsic curvature of $\sigma$ \& $\gamma$
in $R_\gamma$}\label{ss:pf2}

In this subsection, we use $\gamma$ to denote the $N$-geodesic
$\gamma_i$ for $i$ sufficiently large.  We shall prove that the
extrinsic curvature of $\sigma$ and $\gamma$ in the ruled surface
$R_\gamma$ is bounded above by $A^2(1 + \epsilon)f (t))/2$, where
$f=d(\sigma, \gamma)$ and $\epsilon \to 0$ as $\gamma $ approaches
$\sigma$, that is, as $i\to\infty $.

We are going to develop two ways of estimating the length of the
projection to $N$ of a chord $\nu$ of $\sigma$ in $R_\gamma $: an
upper bound for $\ell (\pi_N(\nu ))$, by involving the curvature of
$\sigma$ in $R_\gamma $ to bound the distances $d(\nu(u), N)$; and a
lower bound, by involving that curvature to estimate the distances of
$\pi_N(\nu)$ from its chord $\sigma$ in $N$ (see Figure
\ref{fig:lunes}). The combination of these two bounds will yield the
desired inequality on the curvature of $\sigma$.  The argument is
symmetric, giving the same bound for $\gamma$.

Let $k$ denote the least curvature bound of $\sigma$ in $R_\gamma$,
i. e., for every subarc of $\sigma$ of length $s$ and chord length $r$
in $R_\gamma$ we have $s-r \le \frac{k^2}{24}s^3 + o(s^3)$, and $k$ is
the least such number.

\begin{figure}[h]
\begin{center}
\includegraphics[width= 4.0in]{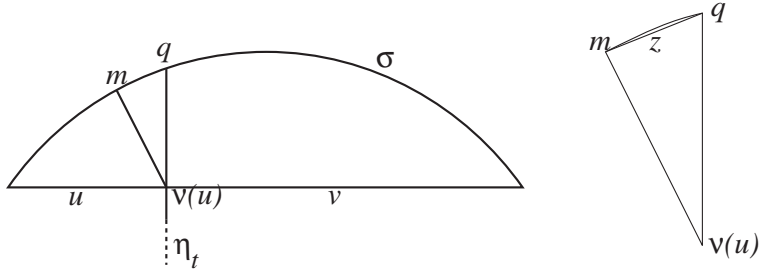}
\end{center}
\caption{Arcs and chords in $R_\gamma$}
\label{fig:curvature}
\end{figure}

By Step 1, we may assume that the rulings $\eta_t$ of $R_\gamma $ make
angles with $\sigma$ and $\gamma$ which differ from $\pi/2$ by at most
$\delta$. Consider a subarc of $\sigma$ with chord $\nu$ of length $r$
in $R_\gamma$. See Figures \ref{fig:lunes} and \ref{fig:curvature}.  A
triangle $\triangle \nu(u) mq$ in $R_\gamma$ formed by a point
$\nu(u)$ on the chord, its nearest point $m$ on $\sigma$, and the
point $q = \sigma(t)$ where the ruling $\eta_t$ through $\nu(u)$ ends
on $\sigma$,
has at least one side of magnitude $d(\nu(u),m) = O(r^2)$, by Lemma
\ref{lem:CBC2}.2.  The other two sides are surely at most $O(r)$;
hence the model area is $O(r^3)$, and the angle sum is bounded by $\pi
+ O(r^3)$.

First we get an upper bound on the length of the projection to $N$ of
$\nu$:

\begin{lem} \label{lem:ubproj}
For a chord $\nu$ of $\sigma$, of chordlength $r$:
\begin{equation} \label{eq:curvature}
\ell(\pi_N\nu) \le r + \frac{r^3}{24}A^2k(1+\delta^2/2 +
O(r^2))f(t)(1+O(f(t)^2)),
\end{equation}
where $f(t)=d(\sigma(t),\gamma(t))$.
\end{lem}

\begin{proof} By Lemma \ref{lem:width}.1,
\begin{equation} \label{eq:width3}
d(\nu(u),m) \le \frac{k}{2}uv + O(r^4).
\end{equation}
(Distances between points of $R_\gamma$ are the distances within
$R_\gamma$.)  Let $\alpha = \angle \nu(u)mq, \, \beta = \angle
\nu(u)qm, \, \theta = \angle m\nu(u)q$, and $z = d(m,q)$. The base
angles formed by the chord $mq$ and the subarc of $\sigma$ with ends
$m, q$ are bounded by $\frac{k}{2}z + O(z^3)$, by Lemma
\ref{lem:CBC2}.1, so that we get lower bounds on $\alpha, \beta$:
\begin{equation} \label{eq:near-right}
\alpha - [\frac{\pi}{2} - \frac{k}{2}z - O(z^3)] \ge 0, \,
\beta - [\frac{\pi}{2} - \delta - \frac{k}{2}z - O(z^3)] \ge 0.
\end{equation}

The angles of the model triangle of $\triangle \nu(u)mq$ in $S_K$
satisfy $\widetilde\alpha \ge \alpha, \, \widetilde\beta \ge \beta, \,
\widetilde\theta \ge \theta$. Add the inequalities
(\ref{eq:near-right}) with $\alpha, \beta$ replaced by
$\widetilde\alpha, \widetilde\beta$ to the obvious inequality
$\widetilde\theta \ge 0$ and 
$\widetilde\alpha + \widetilde\beta + \widetilde\theta \le \pi +
O(r^3)$; then the summands $\pi$ and $-\frac{\pi}{2} -\frac{\pi}{2}$
cancel, leaving $\delta + k +O(z^3) + O(r^3)$ as an upper bound for
the sum of three nonnegative terms. Hence each term must have that
same upper bound, which gives the following.
\begin{equation} \label{eq:tight-right}
\widetilde\alpha \le \frac{\pi}{2} + \delta + \frac{k}{2}z + O(z^3)
+ O(r^3),
\end{equation}
$$\widetilde\beta \le \frac{\pi}{2} + \frac{k}{2}z + O(z^3) + O(r^3),$$
$$\widetilde\theta \le \delta + kz + O(z^3)+ O(r^3).$$ We shall apply
the law of sines to $\triangle \widetilde\nu(u)\widetilde m\widetilde
q$ to estimate $z$ and $d(\nu(u), q) = d(\widetilde\nu(u), \widetilde
q)$.  In doing so we may use the Euclidean version, since the higher
order terms in $\sin(\sqrt{K}z)$ and $\sin(\sqrt{K}d(\nu(u),q))$ can
be absorbed in the error term $O(r^3)$. Thus we have, from
(\ref{eq:width3}),
$$z = \frac{\sin\widetilde\theta}{\sin\widetilde\beta}d(\nu(u), m)
+O(r^3) \le \frac{\delta + kz}{\sin\widetilde\beta}kuv +O(r^3).$$ To
continue we require a lower bound on $\sin\widetilde\beta$.  Note that
the range of $\widetilde\beta$ is bounded above by
(\ref{eq:tight-right}), and below by $\widetilde\beta \ge \beta \ge
\pi/2 -(\delta +\frac{k}{2}z)$ from (\ref{eq:near-right}).  A lower
bound on $\sin\widetilde\beta$ is obtained by checking the endpoints.
On the right, since $z = O(r^2)$, we have $\sin[\frac{\pi}{2} +
\frac{k}{2}z + O(z^3) + O(r^3)] = 1 + O(r^4)$. Hence, we must use the
left, namely, $\sin\widetilde\beta \ge \cos(\delta + \frac{k}{2}z) = 1
- \frac{\delta^2}{2} + O(r^2)$ and $(\sin\widetilde\beta)^{-1} \le 1 +
\frac{\delta^2}{2} + O(r^2)$. This gives:
\begin{equation} \label{eq:short-leg}
z \le kuv(\delta + \frac{\delta^3}{2}) + O(r^4).
\end{equation}
In turn the law of sines also gives $d(\nu(u), q) =
(\sin\widetilde\alpha/\sin\widetilde\beta)d(\nu(u), m)$, and so:
\begin{equation} \label{eq:hypotenuse}
d(\nu(u),q) \le
\frac{1}{1-\delta^2/2}\frac{kuv}{2} + O(r^4) = \frac{k(1 +
\delta^2/2)}{2} uv + O(r^4).
\end{equation}

Since the geodesic in $N$ having $\eta_t$ as chord has curvature $\le
A$ in $M$, we can also apply Lemma \ref{lem:width}.1 to estimate
$d(\nu(u),N)$.  Accordingly, in Lemma \ref{lem:width}.1, replace
$u,r,v$ by $\bar u = d(m,q),\, \bar r = f(t) = d(\sigma(t),
\gamma(t)), \, \bar v = \bar r - \bar u \le f(t)$, obtaining
\begin{eqnarray}\label{eq:distance-N}
d(\nu(u),N) &\le& \frac{A}{2} d(\nu(u),q)f(t)(1 + O(f(t)^2))\\
&\le&
\frac{Ak(1+\delta^2/2 + O(r^2))}{4}uvf(t)(1 + O(f(t)^2)).\nonumber
\end{eqnarray}
Then in Lemma \ref{lem:width}.2 we find that the
original chord of length $r$ has projection to $N$ of length
\begin{eqnarray}
\ell(\pi_N\nu) &\le& r + \int_0^r u(r-u)du \cdot
\frac{A^2k(1+\delta^2/2 + O(r^2))}{4}f(t)(1+O(f(t)^2))\nonumber\\
&=& r + \frac{r^3}{24}A^2k(1+\delta^2/2 + O(r^2))f(t)(1+O(f(t)^2)). \nonumber
\end{eqnarray}
\end{proof}

Now we develop a lower bound for $\ell(\pi_N\nu)$ when $\nu$ is a
chord which gives a sufficiently good approximation for the curvature
bound $k$. For any given positive number $\lambda$ there is a subarc
$\sigma_1$ of length $s_1$ with chord length $r_1$ such that $s_1-r_1
\ge \frac{(k-\lambda)^2}{24}s_1^3$, where $s_1 \to 0$ as $\lambda \to
0$.  Below we assume $\lambda < k/6.$ The length $r_1$ can be attained
from $\sigma_1$ by a process involving first variation: parametrize
$\sigma_1$ by arclength $s \in [0,s_1]$ and let the base angle of the
chord from $\sigma(0)$ to $\sigma(s)$ at $\sigma(s)$ be $\varphi(s) =
\bar \varphi(s) s$. Then the first variation formula gives $r_1$ by
integrating the derivative of the length of those chords, i.e.,
$$r_1 = \int_0^{s_1} \cos \varphi(s) \, ds \ge s_1 -
\sup\frac{\bar\varphi(s)^2}{2} \cdot \frac{s_1^3}{3} + O(s_1^5).$$
Then a point $s_2$ which approximates $\sup\bar\varphi(s)$
sufficiently closely gives us a chord with a relatively large base
angle:
$$\bar\varphi(s_2)^2 \cdot \frac{s_1^3}{6} \ge
\frac{(k-2\lambda)^2}{24} \cdot s_1^3.$$
In this the discrepancy between $\sup\frac{\bar\varphi(s)^2}{2}
\cdot \frac{s_1^3}{3} + O(s_1^5)$ and an upper bound on $s_1 - r_1$ has
been absorbed by changing $\lambda$ to $2\lambda$. Then
\begin{equation}\label{eq:baseang}
\varphi(s_2) \ge \frac{k-2\lambda}{2} s_2.
\end{equation}

\begin{figure}[h]
\begin{center}
\includegraphics[width= 3in]{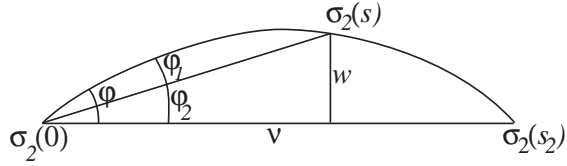}
\end{center}
\caption{Lower bound for $w$}
\label{fig:width}
\end{figure}

Now we turn our attention to the arc $\sigma_2$ of length $s_2$,
parametrized by $s \in [0,s_2]$, with chord $\nu$ of length $r_2$, and
base angle $\varphi = \varphi(s_2)$ satisfying (\ref{eq:baseang}). For
an initial subarc of $\sigma_2$ of length $s$ and base angle
$\varphi_1(s)$, we use the lower bound (\ref{eq:baseang}) on $\varphi$
and the upper bound $\varphi_1(s) \le \frac{k}{2}s + O(s^3)$ of Lemma
\ref{lem:CBC2}.1, to get a lower bound on the distance $w$ from
$\sigma_2(s)$ to $\nu$.  See Figure \ref{fig:width}.
The chord of that subarc is the
hypotenuse (with length $\mu \ge s - \frac{k^2}{24}s^3 + O(s^5)$) of
a right triangle
with angle $\varphi_2(s)$ at $\sigma_2(0)$ and opposite leg of length $w$.
        From the triangle inequality for angles we have $\varphi_2(s) \ge
\frac{k-2\lambda}{2}s_2 - \frac{k}{2}s + O(s^3)$, so that by
side-angle-side comparison we have
\begin{equation}\label{eq:width-lower-bound} w \ge
\mu \cdot \sin\varphi_2(s) \ge \frac{1}{2}\left( ks(s_2-s) -
3\lambda ss_2\right).
\end{equation}

The foot of the leg of length $w$ in the above right triangle projects
to $N$ to give a point at distance $w'$ from $\sigma_2(s)$. The
triangle inequality gives a satisfactory lower bound for $w'$, using
the same bounds on distances to $N$ as in (\ref{eq:distance-N}),
\begin{equation}\label{eq:w-lower-bound-N}
w'(s) \ge w(1-\epsilon) \ge \frac{1-\epsilon}{2} \left( ks(s_2-2) -
3\lambda ss_2\right) = y(s),
\end{equation}
where $\epsilon \to 0$ as $\gamma \to \sigma$.

We let $\nu$ be the chord of $\sigma_2$, and apply Lemma
\ref{lem:shrink} to the curve consisting of $\pi_N\nu$ and geodesic
segment $\sigma_2$ in the space $N$. We may assume that $s_2$ is so
small that the target of $\psi$ can be considered to be the Euclidean
plane, omitting the negligible error terms. If we take the piece of
the parabola $y = y(x)$ for $x \in [0, (1-\sqrt{3\lambda/k})s_2]$ and
its tangent line through $(0,s_2)$ (the point of tangency determines
the interval for the piece), then we have a concave curve shorter than
the length $\psi(\pi_N(\nu))$, which is in turn shorther than
$\pi_N(\nu)$. A lengthy  but straightforward power series
calculation gives a lower bound, for sufficiently small
$s_2$:
$$\ell(\pi_N\nu) \ge s_2 +
(1-\epsilon)^2 \frac{(k-5\lambda)^2}{24}s_2^3.$$

Just as we obtained a lower bound for $\ell(\pi_N\nu)$, we can obtain
a lower bound for the length $s_2$ of $\sigma_2$ by expressing the
distances $w$ in terms of the arclength parameter $r \in [0,r_2]$ of
its chord, which yields the same expression except for higher order
error terms:
$$w \ge \frac{1}{2}\left(kr(r_2-r) - 3\lambda rr_2\right).$$
As before, Lemma \ref{lem:shrink} and an integration to obtain the length
of a concave curve in the Euclidean plane gives
$$s_2 \ge r_2 + \frac {(k-5\lambda)^2}{24}s_2^3.$$
where $s_2 \to 0$ as $\lambda \to 0$ and we have converted $r_2^3$
into $s_2^3$ at the expense of higher order terms. Now we can add
these two inequalities to obtain and cancel the summand $s_2$ to obtain
a lower bound
$$\ell(\pi_N\nu) \ge r_2 + (1-\epsilon)^2 \frac{(k-4\lambda)^2}{24}s_2^3
+ \frac {(k-5\lambda)^2}{24}s_2^3.$$
Now we chain this with the upper bound inequality (\ref{eq:curvature}) with
$r=r_2$, with $r_2^3$ replaced by $s_2^3$. The summands $r_2$ can be
canceled, and in the resulting equation we can divide by $s_2^3$ to leave an
inequality on $k$:
$$(1-\epsilon)^2 \frac{(k-4\lambda)^2}{24} + \frac {(k-5\lambda)^2}{24}
\le \frac{1}{24}A^2k(1+\delta^2/2 + O(r^2))f(t)(1+O(f(t)^2)).$$
In this inequality $4\lambda, 5\lambda$ can be removed by letting
$s\lambda \to 0$.  Then cancelling a factor $k/24$ gives the required
estimate, namely,
$k\le A^2(1 +\epsilon')f (t))/2$, where $\epsilon' \to 0$ as $\gamma
$  approaches $\sigma$.


\subsection{Step 3: $(K+A^2)$-convexity}
\label{ss:pf3}

As in the preceding step, $\gamma$ will denote the $N$-geodesic
$\gamma_i$ for $i$ sufficiently large.  We claim that
$f=d(\sigma,\gamma)$ is almost $K+A^2$-convex, so that the
corresponding ``normal Jacobi field length'' $F$ defined by Equation
(\ref{eq:jaclength}) of Step 1 is in fact ($K+A^2$)-convex.

\begin{figure}[h]
\begin{center}
\includegraphics[width= 1.2in]{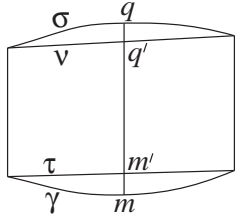}
\end{center}
\caption{Midpoint separations}
\label{fig:strip}
\end{figure}

To estimate the convexity of $f$ we concentrate on the neighborhood of
a single ruling of $R_\gamma$; for convenience, shift the
parametrization so that this ruling is $\eta_0$ and the neighborhood
is the strip between $\eta_{-t}$ and $\eta_t$. The ends of the strip
are bounded by subarcs of $\sigma$ and $\gamma$ and these subarcs will
have geodesic chords in $R_\gamma$, which we denote by $\sigma'$ and
$\gamma'$. See Figure \ref{fig:strip}.
Since $R_\gamma$ has CBA by $K$, the distance between corresponding
points of $\sigma'$ and $\gamma'$ is ($K + \epsilon'$)-convex, where
$\epsilon' \to 0$ as $\gamma \to \sigma$.  Then the distances between
the ends are $f(-t), f(t)$, and with the distance $h$ between the
points $m', q'$ where $\eta_0$ crosses $\gamma', \sigma'$ they must
satisfy the approximate $K + \epsilon'$ midpoint convexity inequality
$$ f(-t) + f(t) - 2h \ge -(K + \epsilon')ht^2 +O(t^3).$$
Now we want to establish the same sort of inequality for $f(-t), f(0),
f(t)$. Since $\eta_0$ is partitioned by $m', q'$ into three segments,
the only difference is the addition of the two end segments, which we
have estimated in Step 2. Letting $u = v = t$ and $m = \gamma(0), q =
\sigma(0)$, we have that $d(m', m)$ and $d(q', q)$ are both bounded by
$A^2(1+\epsilon)f(0)uv/4 + O(t^4) = A^2(1+\epsilon)f(0)t^2/4 +
O(t^4)$. Hence $f(0) = h + d(m', m) + d(q', q) \le h +
A^2(1+\epsilon)f(0)t^2/2 +O(t^4)$. Then
\begin{eqnarray*}
f(-t) + f(t) -2f(0) &\ge& f(-t) + f(t) - 2h - A^2(1+\epsilon)f(0)t^2 +O(t^3)\\
&\ge& -[(K+\epsilon')h + A^2(1+\epsilon)f(0)]t^2 + O(t^3)\\
&=& -(K + A^2 + \epsilon'')f(0)t^2 +O(t^3),
\end{eqnarray*}
where $\epsilon'' \to 0$ as $\gamma \to \sigma$.


\subsection{Step 4: Jacobi field lengths}
\label{ss:pf4}

We have been studying a fan of geodesics in $N$ with lefthand
endpoints fixed at the vertex of a triangle, and righthand endpoints
moving along the opposite side. To conclude that $N$ has CBA by
$K+A^2$, we must establish the local triangle comparison property with
$S_{K+A^2}$. The key is our conclusion in Step 3, according to which
``normal Jacobi field lengths'' $F$ defined as in (\ref{eq:jaclength})
are ($K+A^2$)-convex.  From here, the desired triangle comparison
property may be established by developing the fan into $S_{K+A^2}$,
provided we first establish the appropriate splitting property for
Jacobi field lengths.

Consider reparametrizations of the $\gamma_i$ by $[0,\ell_0]$, setting
$g_i=d(\gamma_i(\ell_i\ell_0^{-1}t),\sigma(t))$ where $\ell_i=
|\gamma_i |$ and $\ell_0=|\sigma|$.  Let $G=\lim_{i\to\infty}
u_i^{-1}g_i$.  We claim
\begin{eqnarray} \label{eq:jacfd1}
G(t)^2 &=& F(t)^2 + t^2 \lim_{i\to\infty} u_i^{-2}(\ell_i -\ell_0)^2\\
&=& F(t)^2 + t^2\cos^2\theta,\nonumber
\end{eqnarray}
where $\theta = \angle pqr$.

The second equality of  (\ref{eq:jacfd1}) is immediate from
the first variation formula.  The first equality is equivalent to
\begin{equation}\label{eq:jacfd2}
g_i(t)^2 - f_i(t)^2 - t^2(\ell_i -\ell_0)^2 = o(u_i^2).
\end{equation}
This expression may be interpreted in terms of the law of cosines for
the planar triangle $\overline{\triangle}$ with sidelengths $f_i(t)$,
$t (\ell_i -\ell_0)$, and $g_i(t)$, namely, the model triangle for
$\triangle = \triangle \sigma(t)\, \gamma_ i(t)\, \gamma_i(t
\ell_0^{-1}\ell_i)$.  Thus the lefthand side of (\ref{eq:jacfd2})
equals
\begin{equation}\label{eq:jacfd3}
f_i(t) t (\ell_i -\ell_0)\cos\overline{\beta}, \end{equation}
where $\overline{\beta}$ is the angle of $\overline{\triangle}$ at
the vertex corresponding to $\gamma_ i(t)$. By triangle comparison,
$\overline{\beta}\ge \beta$, where $\beta$ is the angle of
$\triangle$ at $\gamma_ i(t)$.  By Step 1, $\beta$ is arbitrarily
close to $\pi/2$.

Now apply RMT to the quadrilateral $\Box = \Box \sigma(t)\, \gamma_
i(t)\, \gamma_i(t \ell_0^{-1}\ell_i)\, \sigma(t \ell_0^{-1}\ell_i)$.
The corresponding planar quadrilateral $\overline{\Box}$ has two
adjacent sidelengths agreeing with those of $\overline{\triangle}$,
and their diagonal $\ge g_i(t)$, which is the third sidelength of
$\overline{\triangle}$. Therefore $\widetilde{\beta}\ge \
\overline{\beta}$, where $\widetilde{\beta}$ is the corresponding
angle of $\overline{\Box}$. Since all four angles of $\overline{\Box}$
majorize angles that are, by Step 1, arbitrarily close to $\pi/2$, all
angles of $\overline{\Box}$, including $\widetilde{\beta}$, must be
arbitrarily close to $\pi/2$.  Hence so is $\overline{\beta}$, since
it lies between $\widetilde{\beta}$ and $\beta$.  Thus the expression
(\ref{eq:jacfd3}), when divided by $u_i^2$, approaches $0$, since
$u_i^{-1}f_i(t)\to F(t)$, $u_i^{-1} (\ell_i -\ell_0) \to \cos\varphi$,
and $\cos\overline{\beta} \to 0$.  This completes the verification of
(\ref{eq:jacfd1}).

The triangle comparison property now may be established by development
into $S_{K+A^2}$, as in \cite[p.  709-710]{ABB1}.  The argument
requires that in a neighborhood of each point of $N$, geodesic variations are  Lipschitz if their
endpoint curves are.  This holds because $N$ is known to have some upper 
curvature bound.
Specifically, if $M$ is CAT($K$)  then  $N$ is CAT($k$) for some sufficiently 
large $k$ \cite{L1}.  In this case, we conclude that  $N$ has CAT($K+A^2$) 
neighborhoods of uniform size.


\section{Injectivity radius bounds}
\label{sec:injbd}

\subsection{Injectivity radius }
\label{ss:injdef}

If $M$ is complete and locally compact, the \emph{injectivity radius} $\inj_{M.p}$ is the infimum of sidelengths of \emph{ lunes} (pairs of distinct geodesics with common endpoints) from $p$.
For general inner metric spaces, following a suggestion of Lytchak \cite{L5}, we use a slightly stronger definition of $\inj_{M,p}$, namely:   the supremum of radii $\rho$ with the property that for any $\epsilon >0$ there exists $\delta >0$ such that  any two Lipschitz curves from $p$ to any $q \in B(p,\rho)$, whose lengths are within $\delta$ of $d(p,q)$, are $\epsilon$-close to each 
other. Here we parametrize the curves 
proportionally to arclength by 
$[0,1]$, and use the uniform distance. 
If $M$ is complete, it follows that  balls of radius $< \inj_{M,p}$ possess 
\emph{radial uniqueness}, by which we mean that any point $q$ is joined to $p$ by a unique geodesic and this geodesic varies continuously with $q$. 

It is easily verified that $\inj_{M,p}$ is the infimum of sidelengths of lunes from $p$ in the ultraproduct $M^\omega$.
Here  $M^\omega$ is defined as an 
ultralimit of the sequence $M, M,\ldots$ with fixed basepoint (see \cite{BH}, \cite[\S 2.4]{KL}, \cite{L4}). 
$M^\omega$
is  a complete geodesic space, which contains an isometric copy of $M$, and is isometric to $M$ if $M$ is complete and locally compact.  A sequence of quasi-isometries into a bounded subset of $M$ induces a quasi-isometry into $M^\omega$, and  a complete metric space $M$ is CAT($K$) if and only if 
$M^\omega$ is CAT($K$). 

\subsection{Proof of Theorem  \ref{thm:injrad}}
\label{ss:injpf}

A  \emph{local geodesic} is a locally distance-realizing curve, parametrized proportionally to arclength by $[0,1]$.
In  a complete space of curvature bounded above by $K$, every local 
geodesic $\gamma$ from $p$ of length $<\pi/\sqrt{K}$ has a
neighborhood, in the space of local  geodesics from $p$ in the uniform
topology, on which the righthand endpoint map is injective. The size of this neighborhood is uniformly bounded below in terms of $K$, the length of $\gamma$, and the least size of CAT($K$) neighborhoods of points of $\gamma$.  (See \cite{AB4}, where it is proved that this neighborhood may be taken so that the endpoint map is a homeomorphism onto a neighborhood of the endpoint of $\gamma$.) We say, $M$ has \emph{no conjugate points before $\pi/\sqrt{K}$}.
 
Now suppose $N$ is  a subspace of extrinsic curvature $\le A$ in a CAT($K$) 
space $M$.  Since $N$ has CAT($K+A^2$)  neighborhoods of uniform size, the 
same is true of $N^\omega$ and $(N^\omega)^\omega$.

Suppose, contrary to Theorem  \ref{thm:injrad}, that $\inj_{N,p}< \min\{ \pi/\sqrt{K+A^2},c(A,K)/2 \}$.  There is a sequence of lunes in $N^\omega$ from $p$ with sidelengths approaching $\inj_{N,p}$.  Since
$N^\omega$ has no conjugate points before $\pi/\sqrt{K+A^2}$, the sides are uniformly bounded apart, and so there is a lune from $p$ of sidelength $\inj_{N,p}$  in $(N^\omega)^\omega$. A diagonalization argument gives the existence of such a lune in $N^\omega$.  By Theorem \ref{thm:CBC1}.2, the sides of this lune cannot meet at angle $\pi$, since that would give
 a closed curve in $M$ with extrinsic curvature $\le A$ and length $<c(A,K)$.
Since there is no shorter lune from $p$ in $N^\omega$, the sides also cannot meet at angle $<\pi$. Indeed, in that case they could be deformed to a pair of shorter local geodesics by first
variation;  thus to complete the proof, it is only necessary to verify that a local geodesic $\gamma$ in $N^\omega$ from
$p$ whose length is less than both $\inj_{N,p}$ and $\pi/\sqrt{K+A^2}$ is a geodesic. But if the  maximal minimizing subsegment of $\gamma$ from $p$ had endpoint $\gamma(t), \,\,0<t<1$, there would be geodesics from $p$ to $\gamma(t+\frac{1}{n})$ that were bounded away from this subsegment, again since
$N^\omega$ has no conjugate points before $\pi/\sqrt{K+A^2}$. Then there would be a lune from $p$ in $(N^\omega)^\omega$ of sidelength $< \inj_{N,p}$, and hence by diagonalization a forbidden lune in $N^\omega$.


\subsection{Proof of Corollary \ref{cor:CAT0}}
\label{ss:injcor}

For Part 1, consider a subspace $N$ of extrinsic curvature
$\le A$ in a CAT($0$) space $M$. Since $\pi/\sqrt{A^2}=c(A,0)/2 = \pi/A$, it follows from Theorems \ref{thm:Gausseq} and
\ref{thm:injrad} that $N$ has curvature bounded above by $A^2$ and injectivity radius at least $ \pi/A$.  Therefore balls of smaller radius possess radial uniqueness.  By \cite[Th. 4.3]{AB3}, a ball of radius $\pi/2A$ in $N$ that possesses radial uniqueness is CAT($A^2)$, as required.

In Part 2, $M$ is assumed moreover to be CAT($-A^2$).  Then $N$ has curvature bounded above by $0$ by Theorem
\ref{thm:Gausseq}, and infinite injectivity radius by Theorem
\ref{thm:injrad}, and hence is CAT($0$).  If $N$ were not embedded in
$M$, there would be a nonconstant geodesic in $N$, and hence a curve of 
extrinsic curvature $\le A$ in $M$,  with endpoints 
mapped to a single point of $M$.  Since a complete $A$-curve in the
hyperbolic space of curvature $-A^2$ has infinite length, Theorem
\ref{thm:CBC1}.2 would be contradicted.


\section{Extrinsic curvature of tubular neighborhoods}
\label{sec:sublevels}

Finally we obtain sharp extrinsic curvature bounds for tubular
neighborhoods of convex sets. 

Given the following proposition, the intrinsic curvature bounds in
Example \ref{examp:tub} follow immediately from Theorem
\ref{thm:Gausseq}.  The injectivity radius bounds follow from Theorem
\ref{thm:injrad}, since $\sqrt{1+\tan^2\rho} =\sec\rho$ and
$c(\tan\rho, 1) = 2\pi\cos\rho$.

\begin{prop}\label{prop:tube}
Let $T$ be a $\pi$-totally-convex set in a CAT($1$) space $M$,
and let $N$ be the subset of points at distance $\le \rho$ from $T$,
where $\rho < \pi/2$. Then $N$ has extrinsic curvature $\le \tan\rho$.
\end{prop}

\begin{proof} We reduce to the special case where $T$ is a geodesic
segment in the model space $S_1$. That reduction proceeds as follows.

Let $\nu$ be a geodesic chord of length $<\pi$ in $M$ joining two
points of $N$. Let the endpoints $\nu(0), \nu(r)$ of $\nu$ project to
$m, q \in T$, and $\gamma$ be the geodesic connecting $m$ and $ q$ in $T$,
where $\ell(\gamma) < \pi$. (For the existence and contracting
property of projection to $T$, see \cite[p. 176]{BH}.) By RMT there is
a convex quadrilateral $\Box = \Box\widetilde
m\widetilde\nu(0)\widetilde\nu(r)\widetilde q$ in $S_1$ and a
distance-nonincreasing map $\varphi: \Conv\Box \to M$, preserving the
lengths of the sides of $\Box$. ($Conv$ indicates the convex hull.)
Let $\widetilde\gamma$ be the base of $\Box$, connecting $\widetilde
m, \widetilde q$, and let $\widetilde N$ be the set of points at
distance $\le \rho$ from $\widetilde\gamma$. It is easily seen that
$\widetilde N$ is bounded by two semi-circular arcs of radius $\rho <
\pi/2$ about $\widetilde m, \widetilde q$ with diameters perpendicular
to $\widetilde\gamma$, joined by two circular arcs at distance $\rho$
from the great circle which includes $\widetilde\gamma$. Thus these
joining arcs have curvature $\cot\rho$ and $\widetilde N$ is a subset
of extrinsic curvature $\cot\rho$.  This extrinsic curvature bound is
then inherited by $\widetilde N \cap \Conv\Box$. Let $\widetilde\sigma$
be the geodesic in $\widetilde N$ between $\widetilde\nu(0),
\widetilde\nu(r)$, so $\widetilde\sigma$ has length $\le$ an arc of
curvature $\cot\rho$ and chord-length $r$ and its points are at
distance $\le \rho$ from $\widetilde\gamma$. Using the
distance-nonincreasing property, we have that
$\varphi(\widetilde\sigma)$ has no greater length than
$\widetilde\sigma$ and its points are at distance $\le \rho$ from
$\gamma$, and hence also from $T$, so $\varphi(\widetilde\sigma)$ is
in $N$. This provides the requisite bound on $d_N(\nu(0),\nu(r))$
corresponding to an extrinsic curvature bound $\tan\rho$ for $N$.
\end{proof}

The extrinsic curvature bound in Theorem \ref{prop:tube} is sharp, as
is shown by the example of $T$ being a geodesic segment in $S_1$,
including the possibility of an entire great circle.  The existence of
some curvature bound for a ball of radius $< \pi$ (corresponding to
taking $T$ a ball of radius $ < \pi/2$) was obtained in \cite{L1}.

\begin{rem} Let $f = \sin \sqrt{K} d_T$ for some $\pi/\sqrt{K}$-totally-convex
set $T$ in a CAT$(K)$ space, $K>0$.  It is proved in \cite{AB2} that
$f$ is \emph{$\mathcal{F}K$-convex}, on the subset where $d_T <
\pi/2\sqrt{K}$.  That is, the restriction of $f$ to every geodesic is
$K$-convex. $\mathcal{F}K$-convex functions are basic to the study of spaces 
with curvature bounds (see \cite{AB2, G}). Proposition 
\ref{prop:tube} is a special case of the
following statement, which seems likely to be true, but which we shall
not pursue here: for an $\mathcal{F}K$-convex function $f$ on a
CAT$(K)$ space, $K > 0$, the sublevel set $f^{-1}(-\infty,c]$ for
$c\ge0$ has extrinsic curvature $\le Kc/\inf|Df|$, where the infimum is
taken on the level set $f^{-1}(c)$. (Since the sublevels for $c\le 0$
are convex, the statement is certainly true for $c=0$.) It follows from \cite[Cor 1.9]{L3} that there is some extrinsic curvature bound whenever  the gradient is bounded from below. For $K < 0$,
the corresponding statement would be for nonpositive values of $c$, or
what is the same, for nonnegative superlevel sets of an $\mathcal{F}K$-concave function.

\end{rem}

\section*{Acknowledgments} Our interest in a Gauss Equation in 
Alexandrov spaces of curvature bounded above stems from discussions 
with David Berg and Igor Nikolaev in the mid-90's, for which we thank 
them.  We thank Alexander Lytchak, not only for drawing our attention back by 
his advances in this area, but also for suggesting how to extend Theorem  \ref{thm:injrad} and Corollary \ref{cor:CAT0} to the non-locally-compact case via ultraproducts.


\begin{thebibliography}{ABB}

\bibitem[Ar]{Ar}
S. Alexander, {\em Locally convex hypersurfaces of negatively curved spaces},
Proc. Amer. Math. Soc. {\bf 64} (1977), 321-325.

\bibitem[ABB1]{ABB1} S. B. Alexander, I. D. Berg, R. L. Bishop, {\em
Geometric curvature bounds in Riemannian manifolds with boundary},
Transactions Amer. Math. Soc., {\bf 339} (1993), 703-716.


\bibitem[ABB2]{ABB2} \bysame, {\em
The Riemannian obstacle problem},
Illinois  Math. J., {\bf 31} (1987)2 167-184.


\bibitem[AB1]{AB1}
S. Alexander, R. Bishop, {\em Curvature bounds for warped products of
metric spaces}.  To
appear in Geom. Funct. Anal. (GAFA).
   http://www.math.uiuc.edu/~sba/




\bibitem[AB2]{AB2}
\bysame,
{\em $\mathcal{F}K$--convex functions on metric spaces},  Manuscripta Math.
110 (2003), 115--133.




\bibitem[AB3]{AB3} \bysame, {\em Comparison
Theorems for Curves of
Bounded Geodesic Curvature in Metric Spaces of Curvature Bounded
Above}, Differential Geometry and its Applications, {\bf 6} (1996),
67-86.


\bibitem[AB4]{AB4} \bysame, {\em The
Hadamard-Cartan theorem in locally convex metric spaces},
L'Enseignenment Math., 36(1990), 309-320.



\bibitem[A]{A} A. D. Alexandrov, {\em Ruled surfaces in metric spaces},
Vestnik Leningrad. Univ., 12:5-26, 1957 (Russian).

\bibitem[BH]{BH}
M. Bridson, A. Haefliger,
{\em Metric Spaces of Non-positive Curvature},
Springer-Verlag, Berlin,1999.

\bibitem[BBI]{BBI}
D. Burago, Yu. Burago, S. Ivanov,
{\em A Course in Metric Geometry}, Graduate Studies in Mathematics,
Vol. 33, Amer. Math. Soc., Providence, 2001.



\bibitem[CE]{CE} J. Cheeger,
D. Ebin, {\em Comparison
Theorems in Riemannian Geometry}, North-Holland, Amsterdam, 1975.


\bibitem[Co]{Co} K. Corlette, {\em Immersions with bounded curvature}
Geom. Dedicata {\bf 33} (1990), no. 2, 153--161. MR {\bf 91e}: 53062

\bibitem[EF]{EF}
J. Eells, B. Fuglede, {\em Harmonic Maps between Riemannian Polyhedra}.
With a preface by M. Gromov. Cambridge Tracts in Mathematics, 142.
Cambridge University Press, Cambridge, 2001.

\bibitem[G]{G}
M. Gromov, {\em CAT($\kappa$)-spaces: construction and concentration}
Zap. Nauchn. Sem. S.-Peterburg. Otdel. Mat. Inst. Steklov {\bf
280}(2001), Geom. i Topol. 7, 100-140, 299-300; translation in J. Math.
Sci. (N. Y.) 119 (2004), no. 2, 178--200.


\bibitem[GS]{GS}
M. Gromov, R. Schoen,
{\em Harmonic maps into singular spaces and $p$-adic superrigidity for
lattices in groups of rank one}, Inst. Hautes Etudes Sci.
          Publ. Math. No. 76 (1992), 165--246.


\bibitem[J]{J}
J. Jost, {\em Nonpositive Curvature: Geometric and Analytic Aspects},
Birkhauser, Basel, Boston,1997.

\bibitem[KL]{KL} B. Kleiner, B. Leeb, {\em Rigidity and 
quasi-isometries for symmetric spaces and Euclidean buildings},
Publ. IHES {\bf 86} (1997), 115-197.


\bibitem[L1]{L1} A. Lytchak, {\em Geometry of sets of
positive reach},  Manuscripta Math. {\bf 115} (2004), 199-205.

\bibitem[L2]{L2} \bysame, {\em Almost convex subsets}, to appear in
Geom. Dedicata.

\bibitem[L3]{L3} \bysame, {\em Open map theorem for metric spaces}, to 
appear in St. Petersburg Math. J.

\bibitem[L4]{L4} \bysame, {\em Differentiation in metric spaces}, to 
appear in St. Petersburg Math. J.


\bibitem[L5]{L5} \bysame, {\em Injectivity radius of non-proper spaces}, 
preprint.

\bibitem[M]{M} C. Mese, {\em The curvature of minimal surfaces in
singular spaces}, Comm. Anal. Geom. {\bf 9} (2001), 3-34.

\bibitem[P]{P} A. Petrunin, {\em Metric minimizing surfaces}, Elec.
Res. Announc.  Amer. Math. Soc. {\bf 5} (1998), 47-54.


\bibitem[R]{R}
Yu. G. Reshetnyak,
{\em Nonexpanding maps in a space of curvature no greater than $K$},
Sibirskii Mat. Zh. {\bf 9} (1968), 918-928 (Russian).
English translation: {\em Inextensible mappings in a space of curvature
no greater than $K$}, Siberian Math. Jour. {\bf 9} (1968),683-689.

\bibitem[S]{S} Zh. Shen, {\em A convergence theorem for Riemannian
submanifolds}, Transactions Amer. Math. Soc., {\bf 347} (1995),
1343-1350.

\end{thebibliography}
\end{document}